\documentclass{elsart3-1}

\usepackage{amssymb,amsmath,mathtools}
\usepackage[thmmarks]{ntheorem}
\theoremheaderfont{\bfseries}
\theorembodyfont{\normalfont}
\theoremseparator{ :}
\newtheorem*{proof}{Proof}
\usepackage[english,francais]{babel}

\newtheorem{theorem}{Theorem}[section]
\newtheorem{lemma}[theorem]{Lemma}
\newtheorem{e-proposition}[theorem]{Proposition}

\newtheorem{e-definition}[theorem]{Definition\rm}

\renewcommand{\theequation}{\arabic{equation}}
\def\og{\leavevmode\raise.3ex\hbox{$\scriptscriptstyle\langle\!\langle$~}}
\def\fg{\leavevmode\raise.3ex\hbox{~$\!\scriptscriptstyle\,\rangle\!\rangle$}}

\journal{the Acad\'emie des sciences}
\begin{document}
% place in the next line the header (rubrique) chosen for your article,
% if you know it (you can also have 2, format : Header1/Header2
\centerline{}
\begin{frontmatter}

% Title, authors and addresses

% use the thanksref command within \title, \author or \address for footnotes;
% use the ead command for the email address,
% and the form \ead[url] for the home page:
% \title{Title\thanksref{label1}}
% \thanks[label1]{}
% \author{Name\thanksref{label2}}
% \ead{email address}
% \ead[url]{home page}
% \thanks[label2]{}
% \address{Address\thanksref{label3}}
% \thanks[label3]{}
\selectlanguage{english}
\title{On uniqueness for a rough transport-diffusion equation.}

% use optional labels to link authors explicitly to addresses:
% \author[label1,label2]{}
% \address[label1]{}
% \address[label2]{}
% The [label1] can be suppressed if there is only one address for all authors

\selectlanguage{english}
\author[authorlabel1]{Guillaume L\'evy}
\ead{levy@ljll.math.upmc.fr}

\address[authorlabel1]{Laboratoire Jacques-Louis Lions, Universit\'e Pierre et Marie Curie, Office 15-16 301, Paris}

% If you know the dates of reception, and acceptation you can put them now;
%  idem the name of the person presenting the Note

\medskip
\begin{center}
{\small Received *****; accepted after revision +++++\\
Presented by £££££}
\end{center}

\begin{abstract}
\selectlanguage{english}
% Text of abstract in English
In this Note, we study a transport-diffusion equation with rough coefficients and we prove that solutions are unique in a low-regularity class. 
{\it To cite this article: G. L\'evy, C. R. Acad. Sci. Paris, Ser. I 340 (2005).}

\vskip 0.5\baselineskip

\selectlanguage{francais}
% Text of abstract in French
\noindent{\bf R\'esum\'e} \vskip 0.5\baselineskip \noindent
{\bf Sur l'unicit\'e pour une \'equation de transport-diffusion irr\'eguli\`ere.}
Dans cette Note, nous \'etudions une \'equation de transport-diffusion \`a coefficients irr\'eguliers et nous prouvons l'unicit\'e de sa solution dans une classe de fonctions peu r\'eguli\`eres.
{\it Pour citer cet article~: G. L\'evy, C. R. Acad. Sci. Paris, Ser. I 340 (2005).}

\end{abstract}
\end{frontmatter}

\selectlanguage{english}
% main text
\section{Introduction}
\label{}
In this note, we address the problem of uniqueness for a transport-diffusion equation with rough coefficients. 
Our primary interest and motivation is a uniqueness result for an equation obeyed by the vorticity of a Leray-type solution of the Navier-Stokes equation in the full, three dimensional space.
The main theorem of this note is the following. \hfill\break
\begin{theorem}
 Let $v$ be a divergence free vector field in $L^2(\mathbb{R}_+, \dot{H}^1(\mathbb{R}^3))$ and $a$ a function in $L^2(\mathbb{R}_+ \times \mathbb{R}^3)$.
 Assume that $a$ is a distributional solution of the Cauchy problem
  \begin{equation}  
  (C) \left \{
\begin{array}{c  c}
    \partial_t a + \nabla \cdot  (a v) - \Delta a = 0 \\
    a (0)  = 0,  \\
\end{array}
\right.
 \end{equation}
where the initial condition is understood in the distributional sense. Then $a$ is identically zero on $\mathbb{R}_+ \times \mathbb{R}^3$.
\label{Uniqueness}
\end{theorem}
As a preliminary remark, the assumptions on both $v$ and $a$ entail that $\partial_t a$ belongs to $L^1_{loc}(\mathbb{R}_+, H^{-2}(\mathbb{R}^3))$ and thus, in particular, 
$a$ is also in $\mathcal{C}(\mathbb{R}_+, \mathcal{D}'(\mathbb{R}^3))$.
In Theorem \ref{Uniqueness}, $a$ is to be thought of as a scalar component of the vorticity of $v$, which is in the original problem a Leray solution of the Navier-Stokes equation.
In particular, we only know that $a$ belongs to 
$L^2(\mathbb{R}_+ \times \mathbb{R}^3)$ and $L^{\infty}(\mathbb{R}_+ , \dot{H}^{-1}(\mathbb{R}^3))$, though we will not use the second assumption.
The reader accustomed to three-dimensional fluid mechanics will notice that, comparing the above equation with the actual vorticity equations in $3D$, a term of the type $a\partial_i v$ is missing.
In the original problem, where this Theorem first appeared, we actually rely on a double application of Theorem \ref{Uniqueness}.
For some technical reasons, only the second application of Theorem \ref{Uniqueness} takes in account the abovementioned term.

As opposed to the standard DiPerna-Lions theory, we cannot assume that $a$ is in $L^{\infty}(\mathbb{R}_+ , L^p(\mathbb{R}^3))$ for some $p \geq 1$.
However, our proof does bear a resemblance to the work of DiPerna and Lions; our result may thus be viewed as a generalization of their techniques.
Because of the low regularity of both the vector field $v$ and the scalar field $a$, the use of energy-type estimates seems difficult.
This is the main reason why we rely instead on a duality argument, embodied by the following theorem. \hfill\break
\begin{theorem}
 Given $v$ a divergence free vector field in $L^2(\mathbb{R}_+, \dot{H}^1(\mathbb{R}^3))$ and a smooth $\varphi_0$ in $\mathcal{D}(\mathbb{R}^3)$, there exists a distributional solution of the Cauchy problem 
  \begin{equation}
  (C') \left \{
\begin{array}{c  c}
    \partial_t \varphi - v \cdot \nabla \varphi - \Delta \varphi = 0 \\
    \varphi (0)  = \varphi_0  \\
\end{array}
\right.
 \end{equation}
 with the bounds
 \begin{equation}
  \|\varphi(t)\|_{L^{\infty}(\mathbb{R}^3)} \leq \|\varphi_0\|_{L^{\infty}(\mathbb{R}^3)}
 \end{equation}
 and
  \begin{equation}
   \|\partial_j \varphi(t)\|_{L^2(\mathbb{R}^3)}^2 + \int_0^t \|\nabla \partial_j \varphi(s)\|_{L^2(\mathbb{R}^3)}^2 ds \leq 
 \|\partial_j \varphi_0\|_{L^2(\mathbb{R}^3)}^2 + \|\varphi_0\|_{L^{\infty}(\mathbb{R}^3)}^2 \|\partial_j v\|_{L^2(\mathbb{R}_+ \times \mathbb{R}^3)}^2 
  \end{equation}
for $j = 1,2,3$ and any positive time $t$.
\label{Existence}
\end{theorem}
By reversing the arrow of time, this amounts to build, for any strictly positive $T$, a solution on $[0,T] \times \mathbb{R}^3$ of the Cauchy problem
\begin{equation}
  (-C') \left \{
\begin{array}{c  c}
    - \partial_t \varphi - v \cdot \nabla \varphi - \Delta \varphi = 0 \\
    \varphi (T)  = \varphi_T, \\
\end{array}
\right.
 \end{equation}
where we have set $\varphi_T := \varphi_0$ for the reader's convenience.

\section{Proofs}We begin with the dual existence result. %\hfill\break
\begin{proof}[of Theorem \ref{Existence}.]
Let us choose some mollifying kernel $\rho = \rho(t,x)$ and denote $v^{\delta} := \rho_{\delta} \ast v$, where $\rho_{\delta}(t,x) := \delta^{-4} \rho(\frac{t}{\delta}, \frac{x}{\delta})$. 
Let $(C'_{\delta})$ be the Cauchy problem $(C')$ where we replaced $v$ by $v^{\delta}$.
The existence of a (smooth) solution $\varphi^{\delta}$ to $(C'_{\delta})$ is then easily obtained thanks to, for instance, a Friedrichs method combined with heat kernel estimates.
We now turn to estimates uniform in the regularization parameter $\delta$.
%The first one is simply the usual $L^2$ energy estimate. Multiplying the equation on $\varphi^{\delta}$ by $\varphi^{\delta}$ and integrating in space and time, we get
%\begin{equation}
 %\frac 12 \|\varphi^{\delta}(t)\|_{L^2(\mathbb{R}^3)}^2 + \int_0^t \|\nabla \varphi^{\delta}(s)\|_{L^2(\mathbb{R}^3)}^2 ds = \frac 12 \|\varphi_0\|_{L^2(\mathbb{R}^3)}^2 .
%\end{equation}
The first one is a sequence of energy estimates done in $L^p$ with $p \geq 2$, which yields the maximum principle in the limit. 
Multiplying the equation on $\varphi^{\delta}$ by $\varphi^{\delta} |\varphi^{\delta}|^{p-2}$ and integrating in space and time, we get
\begin{equation}
 \frac 1p \|\varphi^{\delta}(t)\|_{L^p(\mathbb{R}^3)}^p + (p-1) \int_0^t \| \nabla \varphi^{\delta}(s) |\varphi^{\delta}(s)|^{\frac{p-2}{2}} \|_{L^2(\mathbb{R}^3)}^2 ds = \frac 1p \|\varphi_0\|_{L^p(\mathbb{R}^3)}^p .
\end{equation}
Discarding the gradient term, taking $p$-th root in both sides and letting $p$ go to infinity gives
\begin{equation}
 \|\varphi^{\delta}(t)\|_{L^{\infty}(\mathbb{R}^3)} \leq \|\varphi_0\|_{L^{\infty}(\mathbb{R}^3)}.
 \label{BorneInfinie}
\end{equation}
To obtain the last estimate, let us derive for $1 \leq j \leq 3$ the equation satisfied by $\partial_j \varphi^{\delta}$.
We have
\begin{equation}
 \partial_t \partial_j \varphi^{\delta} - v^{\delta} \cdot \nabla \partial_j \varphi^{\delta} - \Delta \partial_j \varphi^{\delta} = \partial_j v^{\delta} \cdot \nabla \varphi^{\delta} .
\end{equation}
Multiplying this new equation by $\partial_j \varphi^{\delta}$ and integrating in space and time gives
\begin{multline}
 \frac 12 \|\partial_j \varphi^{\delta}(t)\|_{L^2(\mathbb{R}^3)}^2 + \int_0^t \|\nabla \partial_j \varphi^{\delta}(s)\|_{L^2(\mathbb{R}^3)}^2 ds = \frac 12 \|\partial_j \varphi_0\|_{L^2(\mathbb{R}^3)}^2 \\ 
 + \int_0^t \int_{\mathbb{R}^3} \partial_j \varphi^{\delta}(s,x) \partial_j v^{\delta}(s,x) \cdot \nabla \varphi^{\delta}(s,x) dx ds .
 \label{Energy}
\end{multline}
Since $v$ is divergence free, the gradient term in the left-hand side does not contribute to Equation (\ref{Energy}).
Denote by $I(t)$ the last integral written above. 
Integrating by parts and recalling that $v$ is divergence free, we have
\begin{align*}
 I(t) & = - \int_0^t \int_{\mathbb{R}^3} \varphi^{\delta}(s,x) \partial_j v^{\delta}(s,x) \cdot \nabla \partial_j \varphi^{\delta}(s,x) dx ds \\
 & \leq \|\varphi_0\|_{L^{\infty}(\mathbb{R}^3)} \int_0^t \|\partial_j v^{\delta}(s)\|_{L^2(\mathbb{R}^3)} \|\nabla \partial_j \varphi^{\delta}(s)\|_{L^2(\mathbb{R}^3)} ds \\
 & \leq \frac 12 \int_0^t \|\nabla \partial_j \varphi^{\delta}(s)\|_{L^2(\mathbb{R}^3)}^2 ds 
 + \frac 12 \|\varphi_0\|_{L^{\infty}(\mathbb{R}^3)}^2 \int_0^t \|\partial_j v^{\delta}(s)\|_{L^2(\mathbb{R}^3)}^2 ds .
\end{align*}
And finally, the energy estimate on $\partial_j \varphi^{\delta}$ reads
\begin{equation}
 \|\partial_j \varphi^{\delta}(t)\|_{L^2(\mathbb{R}^3)}^2 + \int_0^t \|\nabla \partial_j \varphi^{\delta}(s)\|_{L^2(\mathbb{R}^3)}^2 ds \leq 
 \|\partial_j \varphi_0\|_{L^2(\mathbb{R}^3)}^2 + \|\varphi_0\|_{L^{\infty}(\mathbb{R}^3)}^2 \|\partial_j v\|_{L^2(\mathbb{R}_+, \times \mathbb{R}^3)}^2 .
\end{equation}
Thus, the family $(\varphi^{\delta})_{\delta}$ is bounded in $L^{\infty}(\mathbb{R}_+, H^1(\mathbb{R}^3)) \cap L^2(\mathbb{R}_+, \dot{H}^2(\mathbb{R}^3)) \cap L^{\infty}(\mathbb{R}_+ \times \mathbb{R}^3)$.
Up to some extraction, we have the weak convergence of $(\varphi^{\delta})_{\delta}$ in $L^2(\mathbb{R}_+, \dot{H}^2(\mathbb{R}^3))$ and its weak-$\ast$ convergence in 
$L^{\infty}(\mathbb{R}_+, H^1(\mathbb{R}^3)) \cap L^{\infty}(\mathbb{R}_+ \times \mathbb{R}^3)$ to some function $\varphi$.

By interpolation, we also have $\nabla \varphi^{\delta} \rightharpoonup \nabla \varphi$ weakly in $L^4(\mathbb{R}_+, \dot{H}^{\frac 12}(\mathbb{R}^3))$ as $\delta \to 0$.
As a consequence, because $v^{\delta} \to v$ strongly in $L^2(\mathbb{R}_+, \dot{H}^1(\mathbb{R}^3))$ as $\delta \to 0$, the following convergences hold :
$$\Delta \varphi^{\delta} \rightharpoonup \Delta \varphi \text{ in } L^2(\mathbb{R}_+ \times \mathbb{R}^3) ;$$ 
$$v^{\delta} \cdot \nabla \varphi^{\delta} \: , \: \partial_t \varphi^{\delta} \rightharpoonup v \cdot \nabla \varphi \: , \: \partial_t \varphi \text{ in } L^{\frac 43}(\mathbb{R}_+, L^2(\mathbb{R}^3)). $$ 
%$$\partial_t \varphi^{\delta} \rightharpoonup \partial_t \varphi \text{ in } L^{\frac 43}(\mathbb{R}_+, L^2(\mathbb{R}^3)).$$
In particular, such a $\varphi$ is a distributional solution of $(C')$ with the desired regularity.
$\square$
\end{proof}
%\hfill\break
We now state a Lemma which will be useful in the final proof. %\hfill\break
\begin{lemma}
 Let $v$ be a fixed, divergence free vector field in $L^2(\mathbb{R}_+, \dot{H}^1(\mathbb{R}^3))$.
 Let $(\varphi^{\delta})_{\delta}$ be a bounded family in $L^{\infty}(\mathbb{R}_+ \times \mathbb{R}^3)$.
 Let $\rho = \rho(x)$ be some smooth function supported inside the unit ball of $\mathbb{R}^3$ and define $\rho_{\varepsilon} := \varepsilon^{-3} \rho\left(\frac{\cdot}{\varepsilon}\right)$.
 Define the commutator $C^{\varepsilon, \delta}$ by 
 \begin{equation*}
 C^{\varepsilon, \delta}(s,x) :=  v(s,x) \cdot (\nabla \rho_{\varepsilon} \ast \varphi^{\delta}(s))(x) - (\nabla \rho_{\varepsilon} \ast (v(s) \varphi^{\delta}(s)))(x) .
 \end{equation*}
 Then 
 \begin{equation}
   \|C^{\varepsilon, \delta}\|_{L^2(\mathbb{R}_+ \times \mathbb{R}^3)} \leq 
   \|\nabla \rho\|_{L^1(\mathbb{R}^3)} \|\nabla v\|_{L^2(\mathbb{R}_+, \dot{H}^1(\mathbb{R}^3))} \|\varphi^{\delta}\|_{L^{\infty}(\mathbb{R}_+, H^1(\mathbb{R}^3))}.
 \end{equation}
\end{lemma} %\hfill\break
This type of lemma is absolutely not new. 
Actually, it is strongly reminiscent of Lemma II.1 in \cite{DiPerna-Lions} and serves the same purpose.
We are now in position to prove the main theorem of this note. \hfill\break
\begin{proof}[of Theorem \ref{Uniqueness}.]
 Let $\rho = \rho(x)$ be a radial mollifying kernel and define $\rho_{\varepsilon}(x) := \varepsilon^{-3} \rho(\frac{x}{\varepsilon})$.
 Convolving the equation on $a$ by $\rho_{\varepsilon}$ gives, denoting $a_{\varepsilon} := \rho_{\varepsilon} \ast a$,
 \begin{equation}
  (C_{\varepsilon}) \ \ \partial_t a_{\varepsilon} + \nabla \cdot (a_{\varepsilon} v) - \Delta a_{\varepsilon} 
  = \nabla \cdot (a_{\varepsilon} v) - \rho_{\varepsilon} \ast \nabla \cdot (a v) .
 \end{equation}
Notice that even without any smoothing in time, $a_{\varepsilon}$, $\partial_t a_{\varepsilon}$ are in $L^{\infty}(\mathbb{R}_+, \mathcal{C}^{\infty}(\mathbb{R}^3))$ and 
$L^1(\mathbb{R}_+, \mathcal{C}^{\infty}(\mathbb{R}^3))$ respectively, which is enough to make the upcoming computations rigorous.
In what follows, we let $\varphi^{\delta}$ be a solution of the Cauchy problem $(-C'_{\delta})$, with  $(-C'_{\delta})$ being $(-C')$ with $v$ replaced by $v^{\delta}$.
Let us now multiply, for $\delta, \varepsilon > 0$ the equation $(C_{\varepsilon})$ by $\varphi^{\delta}$ and integrate in space and time.
After integrating by parts (which is justified by the high regularity of the terms we have written), we get
\begin{equation*}
 \int_0^T \int_{\mathbb{R}^3} \partial_t a_{\varepsilon}(s,x) \varphi^{\delta}(s,x) dx ds = 
 \langle a_{\varepsilon}(T), \varphi_T \rangle_{\mathcal{D}'(\mathbb{R}^3), \mathcal{D}(\mathbb{R}^3)}  - 
 \int_0^T \int_{\mathbb{R}^3} a_{\varepsilon}(s,x)  \partial_t \varphi^{\delta}(s,x) dx ds
\end{equation*}
and
\begin{multline*}
 \int_0^T \int_{\mathbb{R}^3} \left[ \nabla \cdot (v(s,x) a_{\varepsilon}(s,x)) - \rho_{\varepsilon}(x) \ast \nabla \cdot (v(s,x) a(s,x))\right] \varphi^{\delta}(s,x) dx ds \\
 = \int_0^T \int_{\mathbb{R}^3} a(s,x) C^{\varepsilon, \delta}(s,x) dx ds,
\end{multline*}
where the commutator $C^{\varepsilon, \delta}$ has been defined in the Lemma.
From these two identities, it follows that
\begin{multline*}
 \langle a_{\varepsilon}(T), \varphi_T \rangle_{\mathcal{D}'(\mathbb{R}^3), \mathcal{D}(\mathbb{R}^3)}  
 = \int_0^T \int_{\mathbb{R}^3} a(s,x) C^{\varepsilon, \delta}(s,x) dx ds \\ - 
 \int_0^T \int_{\mathbb{R}^3} a_{\varepsilon}(s,x) \left(- \partial_t \varphi^{\delta}(s,x) - v(s,x) \cdot \nabla \varphi^{\delta}(s,x) - \Delta \varphi^{\delta}(s,x) \right) dx ds.
\end{multline*}
From the Lemma, we know that $(C^{\varepsilon, \delta})_{\varepsilon, \delta}$ is bounded in $L^2(\mathbb{R}^+ \times \mathbb{R}^3)$.
Because $v \cdot \nabla \varphi^{\delta} \to v \cdot \nabla \varphi$ in $L^{\frac 43}(\mathbb{R}^+, L^2(\mathbb{R}^3))$ as $\delta \to 0$, 
the only weak limit point in $L^2(\mathbb{R}_+ \times \mathbb{R}^3)$ of the family $(C^{\varepsilon, \delta})_{\varepsilon, \delta}$ as $\delta \to 0$ is $C^{\varepsilon, 0}$.
Thanks to the smoothness of $a_{\varepsilon}$ for each fixed $\varepsilon$, we can take the limit $\delta \to 0$ in the last equation, which leads to 
\begin{equation}
  \langle a_{\varepsilon}(T), \varphi_T \rangle_{\mathcal{D}'(\mathbb{R}^3), \mathcal{D}(\mathbb{R}^3)} = \int_0^T \int_{\mathbb{R}^3} a(s,x) C^{\varepsilon, 0}(s,x) dx ds .
\end{equation}
Again, the family $(C^{\varepsilon, 0})_{\varepsilon}$ is bounded in $L^2(\mathbb{R}_+ \times \mathbb{R}^3)$ 
and its only limit point as $\varepsilon \to 0$ is $0$, simply because $v \cdot \nabla \varphi_{\varepsilon} - \rho_{\varepsilon} \ast ( v \cdot \nabla \varphi) \to 0$ 
in $L^{\frac 43}(\mathbb{R}^+, L^2(\mathbb{R}^3))$.
Taking the limit $\varepsilon \to 0$, we finally obtain
 \begin{equation}
   \langle a(T), \varphi_T \rangle_{\mathcal{D}'(\mathbb{R}^3), \mathcal{D}(\mathbb{R}^3)} = 0.
 \end{equation}
This being true for any test function $\varphi_T$, $a(T)$ is the zero distribution and finally $a \equiv 0$.
$\square$
\end{proof}
% etc, etc

% The Appendices part is started with the command \appendix;
% appendix sections are then done as normal sections
% \appendix

% \section{}
% \label{}

% The Acknowledgements are an un-numbered section
%\section*{Acknowledgements}
% Acknowledgements text here


\begin{thebibliography}{00}
% please try to use the bibitem system -
% the references should be in alphabetical order of authors' names.
% Articles with a single author first, author will 1 co-author next,
% then author with several co-authors;

% \bibitem{label}
% Text of bibliographic item
\bibitem{AmbrosioBV}
L. Ambrosio, \textsl{Transport equation and Cauchy problem for BV vector fields}, Invent. math. 158, no.2, 227-260 (2004)
\bibitem{DiPerna-Lions} 
R.J. DiPerna and P.-L. Lions, \textsl{Ordinary differential equations, transport theory and Sobolev spaces}, Invent. math. 98, 511-547 (1989)
\bibitem{LebeauFabre}
C. Fabre and G. Lebeau, \textsl{R\'egularit\'e et unicit\'e pour le probl\`eme de Stokes}, Comm. Part. Diff. Eq. 27, no. 3-4, 437-475 (2002)
\bibitem{LeBrisLions}
C. Le Bris and P.-L. Lions, \textsl{Existence and uniqueness of solutions to Fokker-Planck type equations with irregular coefficients}, Comm. Part. Diff. Eq. 33, no. 7-9, 1272-1317 (2008)
\bibitem{Leray} 
J. Leray, \textsl{Sur le mouvement d'un liquide visqueux emplissant l'espace}, Acta Mathematica 63, 193-248 (1934)





\end{thebibliography}
\end{document}